\documentclass[12pt]{amsart}
\usepackage[top=1in, bottom=1in, left=1.5in, right=1.5in]{geometry}
\usepackage{amsmath,amssymb, comment}
\usepackage{caption}
\usepackage{subcaption}
\usepackage{amsthm}
\numberwithin{equation}{section}
\usepackage{amsfonts}
\usepackage{graphicx}
\usepackage{hyperref}
\usepackage{color}
\theoremstyle{plain}
\newtheorem{thm}{Theorem}[section]

\newtheorem{lem}[thm]{Lemma}

\theoremstyle{definition}

\theoremstyle{remark}
\newtheorem{remark}[thm]{Remark}
\author[TY Hou]{Thomas Y Hou}
\address{Applied and Computational Mathematics, California Institute of Technology, Pasadena, CA 91125, USA}
\email{hou@cms.caltech.edu}

\author[P Liu]{Pengfei Liu}
\address{Applied and Computational Mathematics, California Institute of Technology, Pasadena, CA 91125, USA}
\email{plliu@caltech.edu}

\author[F Wang]{Fei Wang}
\address{Department of Mathematics, University of Southern California, Los Angeles, CA 90089, USA}
\email{wang828@usc.edu}

\title[Global Regularity]{Global regularity for a family of 3D models of the axi-symmetric Navier-Stokes equations}
\begin{document}
\maketitle 
\begin{abstract}
	We consider a family of 3D models for the axi-symmetric incompressible Navier-Stokes equations. The models are derived by changing the strength of the convection terms in the axisymmetric Navier-Stokes equations written using a set of transformed variables. We prove the global regularity of the family of models in the case that the strength of convection is slightly stronger than that of the original Navier-Stokes equations, {which demonstrates the potential stabilizing effect of convection. }
\end{abstract}
\section{Introduction and Main Result}
The three-dimensional (3D) Euler and Navier-Stokes equations govern the motion of ideal incompressible fluid in the absence of external forcing:
\begin{equation}
	\label{otarg}
	\mathbf{u}_t+\mathbf{u}\cdot\nabla \mathbf{u}=-\nabla p+\nu\Delta \mathbf{u},\quad \nabla\cdot \mathbf{u}=0.
\end{equation}
Here $\mathbf{u}(x,t): \mathbf{R}^3\times [0, T)\to \mathbf{R}^3$ is the 3D velocity vector of the fluid, and $p(x,t): \mathbf{R}^3\times [0, T)\to \mathbf{R}$ describes the scalar pressure.  The viscous term $\nu \Delta \mathbf{u}$ models the viscous forcing in the fluid. In the case that $\nu=0$, equations \eqref{otarg} are referred to as the Euler equations, and in the case that $\nu>0$, equations \eqref{otarg} are referred to as the Navier-Stokes equations. The divergence-free condition $\nabla \cdot \mathbf{u}=0$ guarantees the incompressibility of the fluid. The Euler and Navier-Stokes equations are among the most fundamental nonlinear partial differential equations (PDEs) in nature yet far from being fully understood. The fundamental question regarding the global regularity of the Euler and Navier-Stokes equations with smooth initial data in the 3D setting remains open, and it is generally viewed as one of the most important open questions in mathematical fluid mechanics; see the surveys \cite{fefferman2006existence, constantin2007euler, gibbon2008three, hou2009blow}.  

The Euler equations have the following scaling-invariance 
\begin{equation}
		\mathbf{u}(x, t)\to\frac{\lambda}{\tau}\mathbf{u}\left(\frac{x}{\lambda}, \frac{t}{\tau}\right), \quad  p(x,t)\to \frac{\lambda^2}{\tau^2}p\left(\frac{x}{\lambda}, \frac{t}{\tau}\right),
	\label{eulerinv}
\end{equation}
and for the Navier-Stokes equations, due to the viscous term, the two-parameter symmetry group in \eqref{eulerinv} is restricted to the following one-parameter group
\begin{equation}
	\mathbf{u}(x, t)\to\frac{1}{\tau^{1/2}}\mathbf{u}\left(\frac{x}{\tau^{1/2}}, \frac{t}{\tau}\right), \quad  p(x,t)\to \frac{1}{\tau}p\left(\frac{x}{\tau^{1/2}}, \frac{t}{\tau}\right). 
	\label{nsinv}
\end{equation}

Smooth solutions to the Euler/Navier-Stokes equations \eqref{otarg} enjoy the following energy identity,
\begin{subequations}
	\label{ape}
\begin{equation} 
	\frac{1}{2}\int |\mathbf{u}(x,t)|^2dx+\nu\int_0^t\int_{\mathbf{R}^3}|\nabla \mathbf{u}(x,s)|^2dxds=\frac{1}{2}\int |\mathbf{u}(x,0)|^2dx,
\end{equation}
which implies the following {\it a priori} estimates for Navier-Stokes:
\begin{equation}
	\label{apeo}
	\frac{1}{2}\int_{\mathbf{R}^3}|\mathbf{u}(x,t)|^2dx,\quad \int_0^t\int_{\mathbf{R^3}}|\nabla \mathbf{u}(x,s)|^2dxds\leq C.
	\end{equation}
The above estimates seem to be the only known coercive {\it a priori} estimates for smooth solutions to the Navier-Stokes \eqref{otarg}. The main difficulty for the global regularity problem of the 3D Navier-Stokes equations lies in the fact that these known {\em a priori} estimates \eqref{apeo} are supercritical with respect to the invariant scaling of the equations \eqref{nsinv}; see \cite{tao2008structure, tao2016finite} for more discussion about this supercritical barrier. For the 3D Euler equations, due to the lack of regularization mechanism (there is no viscosity), to prove the global regularity of the solutions becomes even more challenging.
\end{subequations}

In this work, we consider a family of 3D models for the Navier-Stokes equations with axial symmetry, which is proposed in \cite{hou2017singularity},
\begin{subequations}
\label{hjlequation}
\begin{align}
	u_{1,t}+u^ru_{1,r}+u^zu_{1,z}&=\nu(\partial_r^2+\frac{3}{r}\partial_r+\partial_z^2)u_1+2u_1\phi_{1,z},\\
	\omega_{1,t}+u^r\omega_{1,r}+u^z\omega_{1,z}&=\nu(\partial_r^2+\frac{3}{r}\partial_r+\partial_z^2)\omega_1+(u_1^2)_z,\\
	-(\partial_r^2+\frac{3}{r}\partial_r+\partial_z^2)\phi_1&=\omega_1,
\end{align}
with the Biot-Savart law given by
\begin{equation}
	\label{hjlbs}
	u^r=-\epsilon r\phi_{1,z}, \quad u^z=2\epsilon \phi_1+\epsilon r\phi_{1,r}.
\end{equation}
\end{subequations}
We give the derivation of this model in Section~\ref{derivation} for the sake of completeness.

In \eqref{hjlbs}, the parameter $\epsilon$ characterizes the strength of the convection. The case that $\epsilon=1$ corresponds to the original axi-symmetric Euler/Navier-Stokes equations, and the case $\epsilon=0$ corresponds to the 3D model investigated in \cite{lei2009stabilizing, hou2009partial, hou2014finite}. This family of models was proposed in \cite{hou2017singularity} to study the effect of convection on the depletion of nonlinearity or formation of finite-time singularities. This family of models share several regularity results with the original Euler and Navier-Stokes equations, including an energy identity and two well-known non-blowup criteria. {The numerical results in \cite{hou2017singularity} suggest that the inviscid models with weak convection can develop self-similar singularity and such singularity scenario does not seem to persist as the strength of the convection terms increases, specifically for the original axisymmetric Euler. }

For the family of viscous models with $\epsilon\in [1,2)$, we can obtain a maximum principle for a modified circulation quantity 
$\Gamma^\epsilon=u_1r^{2/\epsilon}$, i.e.,
\begin{equation}
	\label{keyape}
         \|\Gamma^\epsilon\|_{L^\infty}\leq \|\Gamma^\epsilon_0\|_{L^\infty},
\end{equation}
which is subcritical with respect to the invariant scaling for $\epsilon>1$.

The models \eqref{hjlequation} can also be written in a velocity-pressure form as 
\begin{subequations}
	\label{velopress}
 \begin{equation}
	\mathbf{v}_t+\epsilon \mathbf{v}\cdot \nabla \mathbf{v}=-\nabla p+\nu\Delta \mathbf{v}+(2\epsilon-2)\frac{v^\theta v^r \mathbf{e_\theta}}{r},
\end{equation}
where the velocity $\mathbf{v}(x,t)$ is a rescaling of the velocity in model \eqref{hjl}:
\begin{equation}
	\mathbf{v}=\frac{u^r}{\epsilon}\mathbf{e_r}+\frac{u^z}{\epsilon}\mathbf{e_z}+\frac{u^\theta}{\epsilon^{\frac{3}{2}}}\mathbf{e_\theta},\quad u^\theta=ru_1.
\end{equation}
\end{subequations}

Next we state our main result.
\begin{thm}
	\label{mainresult}
	Consider the viscous models \eqref{velopress} with $\epsilon\in (\frac{20}{19},2)$ and data 
	\[ \mathbf{v}(\cdot, 0)\in H^4( \mathbf{R}^3).\]
	Then the solution $\mathbf{v}(x,t)$ is globally regular in time.
\end{thm}

{This result further demonstrates the potential stabilizing effect of the convection terms, which has been demonstrated in the numerical results in \cite{hou2017singularity} that the self-similar singularity of the inviscid models with weak convection does not persist as the strength of the convection terms increases.}

To prove the main result Theorem \ref{mainresult}, we can use $L^p$ estimate for $\omega_1$ and $L^q$ estimate for $u_1$. To control the nonlinear vortex stretching term in the equation of $\omega_1$ using the viscous term, we only need to use the subcritical {\em a priori} estimate \eqref{keyape} and the Hardy inequality, under the condition that $q=2p-\frac{p}{2}\epsilon'$ for some $\epsilon'<\epsilon$. However, for the nonlinear term in the equation of $u_1$, the subcritical {\em a priori} estimate \eqref{keyape} seems insufficient, because it can only control the angular component of the velocity. We use a combination of the supercritical energy estimate \eqref{ape} and the subcritical estimate of $\Gamma^\epsilon$ \eqref{keyape} in the nonlinear term in the equation of $u_1$. To bound the nonlinear term using the viscous term, we need the condition $\epsilon > \frac{20}{19}$ in \eqref{keyyoung}. 

In our proof of the main result in section \ref{mainresult}, we only conduct $L^2$ estimate for $\omega_1$, and using any $L^p$ estimate for $\omega_1$ with $p\in (1, +\infty)$ will lead to the same result under the condition $\epsilon > \frac{20}{19}$. 

The rest of this paper is organized as follows. In section \ref{derivation}, we derive the family of the models that we investigate in this work and list some regularity results for these models. We also give a brief review of recent regularity results for the Navier-Stokes equations with axial symmetry. In section \ref{mainproof}, we prove our main result Theorem \ref{mainresult}.

\section{Derivation of the models and review of the literature}
\label{derivation}
Recently the Euler and Navier-Stokes equations with axial symmetry have attracted a lot of interests. The global regularity problem in this setting remains open although a lot of progress has been made. Let $\mathbf{e_r}$, $\mathbf{e_\theta}$ and $\mathbf{e_z}$ be the standard orthonormal vectors defining the cylindrical coordinates,
\begin{equation*}
	\mathbf{e_r}=(\frac{x_1}{r}, \frac{x_2}{r}, 0)^T,\quad \mathbf{e_\theta}=(\frac{x_2}{r}, -\frac{x_1}{r}, 0)^T,\quad \mathbf{e_z}=(0,0,1)^T,
\end{equation*}
where $r=\sqrt{x_1^2+x_2^2}$ and $z=x_3$. Then the 3D velocity field $\mathbf{u}(x,t)$ is called axi-symmetric if it can be written as 
\[ \mathbf{u}(x,t)=u^r(r, z, t)\mathbf{e_r}+u^\theta(r, z, t)\mathbf{e_\theta}+u^z(r, z, t)\mathbf{e_z},\]
where $u^r$, $u^\theta$ and $u^z$ do not depend on the $\theta$ coordinate.

We denote the axi-symmetric vorticity field $\omega$ as,
\[ \mathbf{\omega}(x,t)=\nabla\times \mathbf{u}(x,t)=\omega^r(r, z, t)\mathbf{e_r}+\omega^\theta(r, z, t)\mathbf{e_\theta}+\omega^z(r, z, t)\mathbf{e_z},\]
and then the Euler and Navier-Stokes equations with axial symmetry can be written using the cylindrical coordinates as 
\begin{subequations}
	\label{oeuler}
	\begin{align}
		u^\theta_t+u^ru^\theta_r+u^zu_z^\theta&=\nu(\Delta-\frac{1}{r^2})u^\theta-\frac{u^ru^\theta}{r},\label{angularv}\\
		\omega^\theta_t+u^r\omega_r^\theta+{u^z}\omega_z^\theta&=\nu(\Delta-\frac{1}{r^2})\omega^\theta+\frac{2}{r}u^\theta u^\theta_z+\frac{u^r\omega^\theta}{r},\label{angularw}\\
		-[\Delta-\frac{1}{r^2}]\phi^\theta&=\omega^\theta,\label{angularbs}
	\end{align}
where the radial and angular velocity fields $u^r(r,z,t)$ and $u^z(r,z,t)$ are recovered from the stream function $\phi^\theta$ based on the Biot-Savart law
\begin{equation}
\label{Biot-Savart-1}
u^r=-\partial_z\phi^\theta,\quad u^z=r^{-1}\partial_r(r\phi^\theta).
\end{equation}
\end{subequations}

Note that the equations for angular velocity \eqref{angularv}, axial vorticity \eqref{angularw}, and the Biot-Savart law \eqref{angularbs}-\eqref{Biot-Savart-1} form a closed system. Equations \eqref{oeuler} have a formal singularity on the axis $r=0$ due to the $\frac{1}{r}$ terms. Using the fact that the angular component $u^\theta(r,z)$, $\omega^\theta(r,z)$ and $\phi^\theta(r,z)$ can all be viewed as odd functions of $r$ \cite{liu2006convergence}, Hou and Li introduced the following transformed variables in \cite{hou2008dynamic},
\begin{equation}
	\label{changevariable}
	u_1=\frac{u^\theta}{r},\quad \omega_1=\frac{\omega^\theta}{r},\quad \phi_1=\frac{\phi^\theta}{r},
\end{equation}
to remove the formal singularity in \eqref{oeuler}. This leads to the following reformulated axi-symmetric Navier-Stokes equations:
\begin{subequations}
\label{nsequation}
\begin{align}
	u_{1,t}+u^ru_{1,r}+u^zu_{1,z}&=\nu(\partial_r^2+\frac{3}{r}\partial_r+\partial_z^2)u_1+2u_1\phi_{1,z},\label{u1eq}\\
	\omega_{1,t}+u^r\omega_{1,r}+u^z\omega_{1,z}&=\nu(\partial_r^2+\frac{3}{r}\partial_r+\partial_z^2)\omega_1+(u_1^2)_z,\label{w1eq}\\
	-[\partial_r^2+\frac{3}{r}\partial_r+\partial_z^2]\phi_1&=\omega_1,\label{bseq}
\end{align}
with the Biot-Savart law given by
\begin{equation}
	\label{axieuler}
	u^r=-r\phi_{1,z}, \quad u^z=2\phi_1+r\phi_{1,r}.
\end{equation}
\end{subequations}

In \cite{hou2017singularity}, a family of 3D models for axi-symmetric Euler and Navier-Stokes equations was proposed by changing the Biot-Savart law \eqref{axieuler}:
\begin{equation}
	\label{hjl}
	u^r=-\epsilon r\phi_{1,z},\quad u^z=2\epsilon \phi_1+\epsilon r\phi_{1,r},
\end{equation}
to study the potential stabilizing effect of the convection terms.

The viscous model \eqref{hjl} enjoys the following scaling-invariance:
\begin{equation}
	\label{vmodels}
	u_1(r, z,t)\to \frac{1}{\tau}u_1\left(\frac{r}{\tau^{\frac{1}{2}}}, \frac{z}{\tau^{\frac{1}{2}}}, \frac{t}{\tau}\right), \quad \omega_1(r, z,t)\to\frac{1}{\tau^{\frac{3}{2}}}\omega_1\left(\frac{r}{\tau^{\frac{1}{2}}},\frac{z}{\tau^{\frac{1}{2}}}, \frac{t}{\tau}\right).
\end{equation}

The modified velocity field \eqref{hjl} is still divergence-free
\[ \nabla\cdot \mathbf{v}=\frac{1}{\epsilon}\left((u^rr)_r+(u^zr)_z\right)=0.\]

It was proved in \cite{hou2017singularity} that the models \eqref{hjl} with $\epsilon\in [0, 2)$ share several regularity results with the original Euler and Navier-Stokes equations, including an energy identity, the conservation of a modified circulation quantity, the BKM non-blowup criterion, and the Prodi-Serrin non-blowup criterion.

	Smooth solutions to the models \eqref{hjl} with $\epsilon\in[0,2)$ enjoy the following energy identity with $u^\theta=ru^1$:
\begin{subequations}
	\label{sae}
\begin{equation}
	\label{energy}
	\begin{split}
		&\frac{1}{2}\frac{d}{dt}\int (u^r)^2+(u^z)^2+\frac{1}{2-\epsilon}(u^\theta)^2 rdrdz\\
		&=-\nu\int |\nabla u^r|^2+|\nabla u^z|^2+\frac{(u^r)^2}{r^2}+\frac{1}{2-\epsilon}\frac{(u^\theta)^2}{r^2} rdrdz.
	\end{split}
\end{equation}

Note that the modified energy functional in \eqref{energy}, 
\begin{equation}
	E_\epsilon=\frac{d}{dt}\int (u^r)^2+(u^z)^2+\frac{1}{2-\epsilon}(u^\theta)^2 rdrdz
\end{equation}
is equivalent to that of the original Euler and Navier-Stokes equations, $E_1$, 
\begin{equation}
	\min(1, \frac{1}{2-\epsilon}) E_1\leq E_\epsilon\leq \max(1, \frac{1}{2-\epsilon}) E_1.
\end{equation}

Based on \eqref{energy}, we have the  following {\em a priori} estimates of the solutions 
\begin{equation}
	\begin{split}
	&\|u^\theta(r, z, t)\|_{L^2},\quad \|u^r(r, z, t)\|_{L^2},\\
	&\int_0^t \|\phi_{1,z}(s)\|_{L^2}^2ds=\int_0^t \frac{1}{\epsilon^2}\|\frac{u^r(s)}{r}\|_{L^2}^2ds\leq C.
\end{split}
\end{equation}
\end{subequations}

We define the modified total circulation $\Gamma^\epsilon$ as 
\begin{subequations}
	\label{tcircul}
\begin{equation}
	\Gamma^\epsilon=u_1r^{2/\epsilon},
\end{equation}
and then $\Gamma^\epsilon$ satisfies the following equation 
\begin{equation}
	\Gamma^\epsilon_t+u^r\Gamma^\epsilon_r+u^z\Gamma^\epsilon_z=\nu \left(\Delta -\frac{2}{r}(\frac{2}{\epsilon}-1)\partial_r+\frac{1}{r^2}\frac{2}{\epsilon}(\frac{2}{\epsilon}-2)\right)\Gamma^\epsilon.
\end{equation}
Then for the inviscid model with $\nu=0$, or the viscous model with $\nu>0$, $\epsilon\geq 1$, we have the following maximum principle
\begin{equation}
	\|\Gamma^\epsilon(r, z, t)\|_{L^\infty}\leq \|\Gamma^\epsilon(r, z, 0)\|_{L^\infty}=\|\Gamma^\epsilon_0\|_{L^\infty}.
\end{equation}
\end{subequations}

For the viscous models with $\epsilon>1$, the quantity $\Gamma^\epsilon$ is indeed subcritical with respect to the invariant scaling of the equations in \eqref{vmodels}, which is the key in our proof of the global regularity result for the models in this paper.

Both the inviscid models and the viscous models enjoy the following BKM type criterion for smooth initial data with decay at infinity. If
\begin{subequations}
	\label{BKM}
\begin{equation}
	\int_0^T\|\nabla\times \mathbf{v}(x, t)\|_{\mathbf{BMO}}dt<+\infty,
\end{equation}
then
\begin{equation}
	\mathbf{v}(x,t)\in L^\infty( H^4(\mathbf{R}^3), [0, T]).
\end{equation}
\end{subequations}

The viscous models also enjoy the Prodi-Serrin type of regularity criterion for smooth initial data with decay at infinity. If 
\begin{subequations}
	\label{prodiserrin}
\begin{equation}
	\mathbf{v}(x,t)\in L^q(L^p(\mathbf{R}^3), (0, T)),\ \frac{3}{p}+\frac{2}{q}=1, \ p\in (3, +\infty],\ q\in [2, +\infty),
\end{equation}
then
\begin{equation}
	\mathbf{v}(x,t)\in L^\infty( H^4(\mathbf{R}^3), [0, T]).
\end{equation}
\end{subequations}

In \cite{hou2017singularity}, convincing numerical evidence is presented to show that  the inviscid models with weak convection could develop stable self-similar singularity on the symmetric axis. The singularity scenario in \cite{hou2017singularity} is different from that at the boundary described in \cite{luo2013potentially-2, hou2015self} in the sense that the center of the singularity region is not stationary but traveling along the symmetric axis. As the strength of the convection terms increases, the self-similar singularity scenario becomes less stable. Such finite-time singularity scenario does not seem to persist for the models with strong convection ($\epsilon \ge \epsilon_0$ for some $\epsilon_0 >0$), specifically the original axi-symmetric Euler equations. These results demonstrate the potential stabilizing effect of the convection terms. In this work, we prove the global regularity of the viscous models when the strength of convection slightly stronger than the original Navier-Stokes, i.e. $\epsilon\in (\frac{20}{19},2)$. The result proved in this work further demonstrates the potential stabilizing effect of convection in axi-symmetric Navier-Stokes equations.

The modified total circulation $\Gamma^\epsilon$ \eqref{tcircul} is subcritical with respect to the scaling \eqref{vmodels} for all $\epsilon>1$. However, the estimate \eqref{tcircul} can only control the angular component of the velocity, and using the technique presented in this work we can only prove the regularity of the models for $\epsilon\in (\frac{20}{19},2)$, not $(1,2)$.

Some important progress has been made regarding the regularity of the axi-symmetric Navier-Stokes equations recently; see \cite{chen2008lower, chen2009lower, koch2009liouville, seregin2009type}, and we mention a few related works below. In \cite{hou2008dynamic}, Hou and Li proposed a 1D model by restricting the equations \eqref{nsequation} to the symmetric axis. Using a cancellation property in the equation for $u_{1,z}$, they proved the global regularity of the 1D model with or without viscosity. In \cite{chen2015regularity}, the cancellation property used in \cite{hou2008dynamic} was further exploited, and several critical regularity criteria concerning only the angular velocity are proved. In particular, the authors of \cite{chen2015regularity} showed that if 
$r^du^\theta\in L^q(L^p(\mathbf{R}^3), (0, T))$ with 
\[ d\in [0, 1),\quad (p, q)\in\{(\frac{3}{1-d},\infty]\times[\frac{2}{1-d},\infty]\},\quad \frac{3}{p}+\frac{2}{q}\leq 1-d,\] then the solutions can be smoothly extended beyond $T$.

In \cite{lei2015criticality}, the global regularity was obtained if $|\Gamma|\leq C|\ln r|^{-2}$, and this result was later improved to $|\Gamma|\leq  C|\ln r|^{-\frac{3}{2}}$ in \cite{wei2016regularity}. The cancellation property in the equation of $u_{1,z}$ is crucial for the results in \cite{chen2015regularity, lei2015criticality, wei2016regularity}. However, for the family of models \eqref{hjl} that we study in this paper, this cancellation is destroyed due to the change of strength in the convection terms in \eqref{hjl}.

\section{Proof of the main result}
\label{mainproof}
In this section we prove the main result Theorem \ref{mainresult}. We need the following Hardy inequality in 1D, see \cite{hardy1952inequalities}.
\begin{lem}\label{lem:hardyinequality}
If $\lambda>1$, $\sigma\neq 1$, $f(r)$ is a nonnegative measurable function, and 
\[ F(r)=\int_0^r f(t)dt,\quad \text{for}\quad \sigma>1,\quad F(r)=-\int_r^\infty f(t)dt,\quad \text{for}\quad \sigma<1,\]
then
\begin{equation}
\label{hardye}
\int_0^\infty r^{-\sigma}F^\lambda dr\leq (\frac{\lambda}{|\sigma-1|})^\lambda \int_0^\infty r^{-\sigma}(rf)^\lambda dr.
\end{equation} 
\end{lem}

We also need the following elliptic estimates \cite{hou2008global, miao2013global, lei2015axially, wei2016regularity}
\begin{lem}
	For axi-symmetric smooth functions $\phi_1(r,z)$ and $\omega_1(r,z)$ in $\mathbf{R}^3$, which satisfy the elliptic equation
	\[-\Delta \phi_1-\frac{2}{r}\partial_r\phi_1=\omega_1,\] 
	we have the following estimates
	\begin{equation}
		\label{ellipticestimate}
		\|\nabla^2 \phi_1\|_{L^2}\leq C\|\omega_1\|_{L_2},\quad \|\nabla^2 \phi_{1,z}\|_{L^2}\leq C\|\nabla \omega_1\|_{L^2}.
	\end{equation}
\end{lem}
\begin{subequations}
\begin{lem}\label{lemm}
	For smooth solution of the model \eqref{hjl}, $u_1$ and $\epsilon'\in(1, \epsilon)$, 
	\begin{equation}
		\label{hardy}
		\int |u_1(r, z)|^{\epsilon'} f(r)^2rdrdz\leq C_1(r)\int |\partial_rf|^2rdrdz+Cr_1^{-\frac{2\epsilon'}{\epsilon}}\int_{r\geq r_1}f^2rdrdz,
	\end{equation}
	with
	\begin{equation}
		\label{c1coeff}
		C_1(r_1)=C\|\Gamma_0^\epsilon\|^{\epsilon'}_{L^\infty}r_1^{2-\frac{2\epsilon'}{\epsilon}}\left(\frac{\epsilon}{\epsilon-\epsilon'}\right)^2,\quad \lim_{r_1\to 0^+}C_1(r_1)=0 .
	\end{equation}
\end{lem}
\end{subequations}
\begin{proof} Let $\psi(r)$ be a radial cutoff function such that
	\begin{equation}
		\label{cutoff}
		\psi(r)\in C^\infty (R),\ \psi(r)=\begin{cases}
			1,\quad r\leq 1\\
			0,\quad r\geq 2
		\end{cases},\ 0\leq \psi(r)\leq 1,\ |\psi_r(r)|\leq 2.
	\end{equation}
		Denote $\psi_{r_1}(r)$ as $\psi(\frac{r}{r_1})$, then we have 
		\begin{equation}
			\label{split1}
			\begin{split}
				&\int |u_1| ^{\epsilon'}f(r)^2rdrdz=\int |u_1|^{\epsilon'}\left(f(r)\psi_{r_1}(r)+f(r)(1-\psi_{r_1}(r))\right)^2 rdrdz\\
				&\leq 2\int_{r\leq 2r_1} |u_1|^{\epsilon'}f(r)^2|\psi_{r_1}(r)|^2rdrdz+2\int_{r\geq r_1}|u_1|^{\epsilon'}f(r)^2(1-\psi_{r_1}(r))^2rdrdz.
		\end{split}
	\end{equation}

Using the maximum principle \eqref{tcircul}, we have 
	\begin{equation}
		\label{bu1}
		|u_1(r,z,t)|\leq \|\Gamma^\epsilon\|_{L^\infty}r^{-\frac{2}{\epsilon}}\leq \|\Gamma_0^\epsilon\|_{L^\infty}r^{-\frac{2}{\epsilon}}.
	\end{equation}

	Putting \eqref{bu1} in the first term on the RHS of \eqref{split1}, and using the Hardy inequality \eqref{hardye}, we get
	\begin{equation}
		\label{l1p}
		\begin{split}
			&\quad 	\int_{r\leq 2r_1} |u_1|^{\epsilon'}f(r)^2|\psi_{r_1}(r)|^2rdrdz \\
			&\leq \int_{r\leq 2r_1}\|\Gamma_0^\epsilon\|_{L^\infty}^{\epsilon'}r^{-\frac{2\epsilon'}{\epsilon}}|f(r)\psi_{r_1}(r)|^2rdrdz\\
			&\leq C\|\Gamma_0^\epsilon\|^{\epsilon'}_{L^\infty}\left(\frac{\epsilon}{\epsilon-\epsilon'}\right)^2\int_{r\leq 2r_1} r^{2-\frac{2\epsilon'}{\epsilon}}|\partial_r(f(r)\psi_{r_1}(r))|^2rdrdz\\
			&\leq C\|\Gamma_0^\epsilon\|^{\epsilon'}_{L^\infty}\left(\frac{\epsilon}{\epsilon-\epsilon'}\right)^2\left[\int_{r\leq 2r_1}r^{2-\frac{2\epsilon'}{\epsilon}}|\partial_r f|^2rdrdz+\int_{r\geq r_1}r^{2-\frac{2\epsilon'}{\epsilon}}|f|^2|\partial_r\psi_{r_1}|^2rdrdz\right]\\
			&\leq C\|\Gamma_0^\epsilon\|^{\epsilon'}_{L^\infty}\left(\frac{\epsilon}{\epsilon-\epsilon'}\right)^2r_1^{2-\frac{2\epsilon'}{\epsilon}}\|\partial_rf\|^2_{L^2}+C\|\Gamma_0^\epsilon\|^{\epsilon'}_{L^\infty}\left(\frac{\epsilon}{\epsilon-\epsilon'}\right)^2r_1^{-\frac{2\epsilon'}{\epsilon}}\int_{r\geq r_1}f^2rdrdz.
		\end{split}
	\end{equation}

	For the second term in \eqref{split1}, using the estimate \eqref{bu1}, we have 
	\begin{equation}
		\label{l2p}
		\int_{r\geq r_1}|u_1|^{\epsilon'}f^2(r)(1-\psi_{r_1}(r))^2rdrdz\leq C\|\Gamma_0^\epsilon\|^{\epsilon'}_{L^\infty}r_1^{-\frac{2\epsilon'}{\epsilon}}\int_{r\geq r_1}f^2rdrdz.
	\end{equation}

	Adding up estimates \eqref{l1p} and \eqref{l2p}, we prove \eqref{hardye}.  
	\end{proof}

Next we give the proof for the Theorem \ref{mainresult}. Without loss of generality, we assume that $\nu=1$ in our proof.
\begin{proof}
	We denote $\epsilon'=\frac{20}{19}<\epsilon$ and consider the following two quantities: 
	\[\int |\omega_1|^2rdrdz,\quad \int |u_1|^{q}rdrdz,\quad q=4-\epsilon'.\]

	Multiplying the equation of $\omega_1$ \eqref{w1eq} by $\omega_1$, we get
	\begin{equation}
		\label{weqe}
	\begin{split}
		&\frac{d}{dt}\frac{1}{2}\int \omega_1^2rdrdz+\frac{1}{2}\int u^r (\omega_1^2)_r+u^z(\omega_1^2)_zrdrdz\\
		&=\int 2u_1u_{1,z}\omega_1rdrdz+\int (\Delta \omega_1+\frac{2}{r}\omega_{1,r})\omega_1rdrdz.
	\end{split}
	\end{equation}

	Using integration by part, we can show that the convection terms vanish due to the incompressibility condition $(u^rr)_r+(u^zr)_z=0$:
	\[\int u^r (\omega_1^2)_r+u^z(\omega_1^2)_zrdrdz=-\int \omega_1^2((u^rr)_r+(u^zr)_z)drdz=0.\]

	For the viscous term on the RHS of \eqref{weqe}, we have
	\begin{equation}
		\label{laplace}
		\int\Delta \omega_1 \omega_1rdrdz=-\int|\nabla \omega_1|^2rdrdz.
	\end{equation}

	Next we treat the first order derivative term on the RHS of \eqref{weqe} as
	\begin{equation}
		\label{firstorder}
		\int \frac{2}{r}\omega_{1,r}\omega_1rdrdz=\int (\omega_1^2)_rdrdz=-\int \omega_1(0, z, t)^2dz\leq 0.
	\end{equation}

Using integration by part and Young's inequality leads to
\begin{equation}
	\label{nonlinearw}
	\begin{split}
	& \left|\int 2u_1u_{1,z}\omega_1rdrdz\right|=\left|\int u_1^2\omega_{1,z}rdrdz\right|\\
	& \leq \frac{1}{2}\int u_1^4 rdrdz+\frac{1}{2}\int |\nabla\omega_1|^2rdrdz.
\end{split}
\end{equation}

For the first term on RHS of \eqref{nonlinearw}, using Lemma \ref{lemm} and $q=4-\epsilon'$,  we get
\begin{equation}
\label{eu4}
\begin{split}
&\frac{1}{2}\int u_1^4rdrdz=\frac{1}{2} \int |u_1|^{\epsilon'}|u_1|^{q}rdrdz\\
&\leq C(r_1)\int |\nabla(|u_1|^{\frac{q}{2}})|^2rdrdz+C\int |u_1|^{q}rdrdz.
\end{split}
\end{equation}
	
Adding up estimates \eqref{laplace} and \eqref{eu4} in \eqref{weqe}, we have 
\begin{equation}
\label{few}
\begin{split}
	 &\frac{d}{dt}\int \omega_1^2rdrdz+\int |\nabla \omega_1|^2rdrdz\\
	 &\leq C(r_1) \int |\nabla(|u_1|^{\frac{q}{2}})|^2rdrdz+C\int |u_1|^{q} rdrdz.
\end{split}
\end{equation}

Next we consider the equation of $u_1$, \eqref{u1eq} and multiply both sides by $|u_1|^{q-2}u_1$ to obtain
\begin{equation}
	\label{u1q}
	\begin{split}
	&\frac{d}{dt}\frac{1}{q}\int |u_1|^{q}rdrdz+\frac{1}{q}\int u^r|u_1|^{q}_r+u^z|u_1|^{q}_zrdrdz\\
	=&\int 2|u_1|^{q}\phi_{1,z}rdrdz+\int (\Delta u_1+\frac{2}{r}u_{1,r})|u_1|^{q-2}u_1rdrdz.
\end{split}
\end{equation}

	Again the convection terms vanish due to incompressibility. For the diffusion term, using estimates similar to those for the $\omega_1$ equation, we arrive at
	\begin{equation}
		\begin{split}
		&\int (\Delta u_1+\frac{2}{r}u_{1,r})|u_1|^{q-2}u_1rdrdz\\
		&=-\frac{4(q-1)}{q^2}\int |\nabla (|u_1|^{\frac{q}{2}})|^2rdrdz-\frac{q}{2}\int |u_1(0, z)|^qdz\\
		&\leq -\frac{4(q-1)}{q^2}\int |\nabla (|u_1|^{\frac{q}{2}})|^2rdrdz.
	\end{split}
	\end{equation}

	Next we decompose the nonlinear term on the RHS of \eqref{u1q} into two parts 
	\begin{equation}
	\label{cut2}
	\int 2|u_1|^{q}\phi_{1,z}rdrdz=2\int |u_1|^q\phi_{1,z}\psi(r)+|u_1|^{q}\phi_{1,z}(1-\psi(r))rdrdz,
	\end{equation}
	where $\psi(r)$ is the cut-off function defined in \eqref{cutoff} that satisfies $\psi(r)=1$ for $r\leq 1$.

	For the second term on the RHS of \eqref{cut2},  Young's inequality implies,
	\begin{equation}
		\label{fareu}
		\begin{split}
&	\int_{r\geq 1} |u_1|^q\phi_{1,z}(1-\psi(r))rdrdz\\
&	\leq \int_{r\geq 1} |u_1|^{2q}+|\phi_{1,z}|^2rdrdz\leq \int_{r\geq 1} |ru_1|^2 r^{-2}|u_1|^{2q-2}+|\phi_{1,z}|^2rdrdz\\
& \leq \|\Gamma_0^\epsilon\|^{2q-2}_{L^\infty}\int_{r\geq 1} |ru_1|^2rdrdz+\frac{1}{\epsilon^2}\int_{r\geq 1} |\frac{u^r}{r}|^2rdrdz\leq C,
\end{split}
	\end{equation}
	where we have used the {\em a priori} estimates \eqref{sae} in the last step. 

	As for the first term on the RHS of \eqref{cut2}, we denote 
\begin{equation}
	\label{defg}
	g(r,z)=\phi_{1,z}(r,z)\psi(r)
\end{equation}
and have 
	\begin{equation}
		\label{facu}
		\left|\int 2|u_1|^{q}g(r,z) rdrdz\right|	\leq 2\int |u_1|^\alpha |u_1|^\beta |g(r,z)|rdrdz,
	\end{equation}
	where the exponents $\alpha$, $\beta$ are
	\begin{equation}
	\alpha=\frac{16-8\epsilon'}{4-\epsilon'},\quad \beta=\frac{(\epsilon')^2}{4-\epsilon'},\quad \alpha+\beta=q=4-\epsilon'.
	\end{equation}
	
	Then applying Young's inequality with 
	\[ L_1=\frac{4-\epsilon'}{4-2\epsilon'},\quad L_2=\frac{4-\epsilon'}{\epsilon'},\quad \frac{1}{L_1}+\frac{1}{L_2}=1,\]
	we obtain 		
\begin{equation}
	\begin{split}
	\label{young1}
	\int 2|u_1|^q|g|rdrdz&	\leq \frac{1}{2}\int |u_1|^{L_1\alpha} rdrdz+C\int |u_1|^{L_2\beta} |g|^{L_2}rdrdz\\
	&=\frac{1}{2}\int |u_1|^4rdrdz+C\int |u_1|^{\epsilon'} |g(r,z)|^{\frac{4-\epsilon'}{\epsilon'}}rdrdz.
\end{split}
\end{equation}
	
	The first term in \eqref{young1} is treated as in the estimate \eqref{eu4}. For the second term on the RHS of \eqref{young1}, using Lemma \ref{lemm} with $r_1=2$ and the fact that $g(r,z)=0$ for $r\geq 2$, we obtain   
	\begin{equation}
		\label{allin}
		\begin{split}
			&\int |u_1|^{\epsilon'}|g(r,z)|^{\frac{4-\epsilon'}{\epsilon'}}rdrdz\leq C\int |\partial_r(|g(r,z)|^{\frac{4-\epsilon'}{2\epsilon'}})|^2rdrdz\\
			& =C \int |g_{r}(r,z)|^2|g(r,z)|^{\frac{4-3\epsilon'}{\epsilon'}}rdrdz \leq C\|g(r,z)\|_{L^\infty}^{\frac{4-3\epsilon'}{\epsilon'} }\|\nabla g(r,z)\|_{L^2}^2 .
		\end{split}
		\end{equation}
\begin{subequations}
\label{interpf}
Using the following interpolation inequality,
\begin{align}
	\|\nabla g(r,z)\|_{L^2}&\leq C\|g(r,z)\|_{L^2}^{\frac{1}{2}}\|\nabla^2g(r,z)\|_{L^2}^{\frac{1}{2}},\\
	\|g(r,z)\|_{L^\infty}&\leq C\|g(r,z)\|_{L^2}^{\frac{1}{4}}\|\nabla^2g(r,z)\|_{L^2}^{\frac{3}{4}},
\end{align}
we have 
\begin{equation}
\|\nabla g(r,z)\|_{L^2}^{\frac{11\epsilon'-4}{6\epsilon'}}\| g(r,z)\|_{L^\infty}^{\frac{4-3\epsilon'}{\epsilon'}}\leq C\|g(r,z)\|_{L^2}^{\frac{\epsilon'+4}{6\epsilon'}}\|\nabla^2 g(r,z)\|_{L^2}^{\frac{8-4\epsilon'}{3\epsilon'}}.
\end{equation}
\end{subequations}
Using \eqref{interpf} in \eqref{allin}, we have 
\begin{equation} 
	\label{fallin}
	\begin{split}
		&\int |u_1|^{\epsilon'}|g(r,z)|^{\frac{4-\epsilon'}{\epsilon'}}rdrdz\leq  C\|g(r,z)\|_{L^\infty}^{\frac{4-3\epsilon'}{\epsilon'} }\|\nabla g(r,z)\|_{L^2}^2 \\
		&= C\left[\|g(r,z)\|_{L^\infty}^{\frac{4-3\epsilon'}{\epsilon'}}\|\nabla g(r,z)\|_{L^2}^{\frac{11\epsilon'-4}{6\epsilon'}}\right] \|\nabla g(r,z)\|_{L^2}^{\frac{4+\epsilon'}{6\epsilon'}}\\
		&\leq C\|g(r,z)\|_{L^2}^{\frac{4+\epsilon'}{6\epsilon'}}\|\nabla g(r,z)\|_{L^2}^{\frac{4+\epsilon'}{6\epsilon'}}\|\nabla^2 g(r,z)\|_{L^2}^{\frac{8-4\epsilon'}{3\epsilon'}}.
	\end{split}
\end{equation}

In deriving the above estimate, we have used the interpolation inequality \eqref{interpf} such that the exponents for $\|g(r,z)\|_{L^2}$ and $\|\nabla g(r,z)\|_{L^2}$ are the same in the RHS of \eqref{fallin}.

Since $\epsilon'=\frac{20}{19}$, using Young's inequality with 
\begin{equation}
	\label{keyyoung}
	L_1=\frac{12\epsilon'}{4+\epsilon'}, \quad L_2=\frac{6\epsilon'}{8-4\epsilon'},\quad \frac{1}{L_1}+\frac{1}{L_2}=1,
\end{equation}
	in \eqref{fallin}, we have 
\begin{equation}
	\label{young2}
	\begin{split}
		&\int |u_1|^{\epsilon'}|g(r,z)|^{\frac{4-\epsilon'}{\epsilon'}}rdrdz\\
		&\leq C(\delta) \left(\|g(r,z)\|_{L^2}^{\frac{4+\epsilon'}{6\epsilon'}}\|\nabla g(r,z)\|_{L^2}^{\frac{4+\epsilon'}{6\epsilon'}}\right)^{L_1}+\delta \left( \|\nabla^2 g(r,z)\|_{L^2}^{\frac{8-4\epsilon'}{3\epsilon'}}\right)^{L_2}\\
		&=C(\delta)\|g(r,z)\|_{L^2}^2\|\nabla g(r,z)\|_{L^2}^2+\delta \|\nabla^2 g(r,z)\|_{L^2}^2.
	\end{split}
\end{equation}

	Since $g(r,z)=\phi_{1,z}\psi(r)$ and $\psi(r)$ is constant for $r\leq 1$, we have 
	\begin{subequations}
		\label{expand}
		\begin{equation}
	\|g(r,z)\|_{L^2}^2\leq C\|\phi_{1,z}\|^2_{L^2},
\end{equation}
	\begin{equation}
		\|\nabla g(r,z)\|_{L^2}^2\leq C\|\nabla \phi_{1,z}\|_{L^2}^2+C\|1_{r\geq 1}\phi_{1,z}\|_{L^2}^2,
\end{equation}
\begin{equation} \|\nabla^2 g(r,z)\|_{L^2}^2\leq C\|\nabla^2\phi_{1,z}\|^2+C\|1_{r\geq 1}\nabla \phi_{1,z}\|_{L^2}^2+C\|1_{r\geq 1}\phi_{1,z}\|_{L^2}^2.
\end{equation}
\end{subequations}
	
By the {\em a priori} estimate \eqref{sae}, we have
\begin{equation}
	\label{lsae}
	\|1_{r\geq 1}\phi_{1,z}\|_{L^2}^2=\|1_{r\geq 1}\frac{u^r}{r}\|_{L^2}^2\leq \|u^r\|_{L^2}^2\leq C.
\end{equation}
In view of \eqref{young2}, we get from \eqref{expand} and \eqref{lsae}
\begin{equation}
	\label{fsob}
	\begin{split}
		&\int |u_1|^{\epsilon'}|g(r,z)|^{\frac{4-\epsilon'}{\epsilon'}}rdrdz \leq C\delta \|\nabla^2\phi_{1,z}\|_{L^2}^2+C\|\nabla \phi_{1,z}\|_{L^2}^2\|\phi_{1,z}\|_{L^2}^2\\
		&+C\|\nabla \phi_{1,z}\|_{L^2}^2+C\|\phi_{1,z}\|_{L^2}+C.
	\end{split}
\end{equation}
Employing the elliptic estimate \eqref{ellipticestimate} in \eqref{fsob}, we deduce 
\begin{equation}
	\label{ffsob}
	\begin{split}
		&\int |u_1|^{\epsilon'}|g(r,z)|^{\frac{4-\epsilon'}{\epsilon'}}rdrdz \leq C\delta \|\omega_{1,z}\|_{L^2}^2+C\|\omega_{1}\|_{L^2}^2\|\phi_{1,z}\|_{L^2}^2\\
		&+C\|\omega_{1}\|_{L^2}^2+C\|\phi_{1,z}\|_{L^2}+C.
	\end{split}
\end{equation}
Putting the estimates \eqref{eu4}, \eqref{fareu}, \eqref{ffsob} in \eqref{u1q}, we get 
\begin{equation}
	\begin{split}
	&\frac{d}{dt}\int |u_1|^qrdrdz+ C(1-C_1(r_1))\|\nabla(|u_1|^{\frac{q}{2}})\|^2_{L^2}-C\delta\|\nabla \omega_{1}\|_{L^2}^2\\
	& \leq 
	C\|\phi_{1,z}\|_{L^2}^2\|\omega_1\|_{L^2}^2+C\|\phi_{1,z}\|_{L^2}^2+C\|\omega_1\|_{L^2}^2+C.
\end{split}
	\label{feu}
\end{equation}

At last, choosing $r_1$, $\delta$ small enough, and adding up \eqref{few} with \eqref{feu} give
\begin{equation}
	\label{finale}
	\begin{split}
		&\frac{d}{dt}\int |\omega_1|^2+|u_1|^qrdrdz+C\int |\nabla \omega_1|^2+|\nabla(|u_1|^{\frac{q}{2}})|^2rdrdz\\
	&\leq C(\|\phi_{1,z}\|_{L^2}^2+1) (\int |\omega_1|^2+|u_1|^q rdrdz)+C\|\phi_{1,z}\|^2+C,
\end{split}
\end{equation}
which together with the {\em a priori} estimate \eqref{sae} implies 
\begin{equation}
	\int \omega_1^2(T)+|u_1(T)|^{q}rdrdz+\int_0^T\int |\nabla \omega_1|^2+|\nabla(|u_1|^{\frac{q}{2}})|^2rdrdzds\leq C
	\label{finalestimateuw}
\end{equation}
where the constant $C$ may depend on the initial data and $T$.

To prove the global regularity of the solutions, we consider 
\begin{equation}
	\label{ev4}
	\|\mathbf{v}(x)\|_{L^4}\leq C\|u^\theta\|_{L^4} +C\|u^r\|_{L^4}+C\|u^z\|_{L^4} .
\end{equation}
For the $\|u^\theta\|_{L^4}$ term in \eqref{ev4}, using estimates \eqref{sae}, \eqref{tcircul}, and \eqref{finalestimateuw}, we obtain
\begin{equation}
	\label{utheta4e}
	\begin{split}
		\int u_1^4r^4rdrdz&=\int_{r\leq 1} |u_1|^{q} |u_1|^{\epsilon'}r^4 rdrdz+\int_{r\geq 1} (u^\theta)^2 \frac{|\Gamma^\epsilon|^2}{r^{\frac{4}{\epsilon}-2}}rdrdz\\
		&\leq \|\Gamma_0^\epsilon\|_{L^\infty}^{\epsilon'}\int |u_1|^qrdrdz+\|\Gamma^\epsilon_0\|_{L^\infty}^2\int (u^\theta)^2rdrdz\leq C.
	\end{split}
\end{equation}

Then we consider the equation for $\omega^\theta=r\omega_1$, which is 
\begin{equation}
	\omega^\theta_t+u^r\omega^\theta_r+u^z\omega^\theta_z=\frac{u^r}{r}\omega^\theta+\frac{(u^\theta)^2_z}{r}+(\Delta-\frac{1}{r^2})\omega^\theta.
	\label{eqwt}
\end{equation}
Multiplying both sides of \eqref{eqwt} by $\omega^\theta$ and integrating, we get 
\begin{equation}
\label{eqwtheta}
\begin{split}
	&\frac{1}{2}\frac{d}{dt}\int (\omega^\theta)^2rdrdz+\frac{1}{2}\int u^r(\omega^\theta)^2_r+u_z(\omega^\theta)^2_zrdrdz\\
	& =\int \frac{u^r}{r}(\omega^\theta)^2rdrdz+\int \frac{(u^\theta)^2_z \omega^\theta}{r}rdrdz-\int |\nabla \omega^\theta|^2+\frac{(\omega^\theta)^2}{r^2}rdrdz.
\end{split}
\end{equation}

The convection terms vanish due to the incompressibility condition, and for the first nonlinear term in \eqref{eqwtheta}, we have 
\begin{equation}
	\label{firstomegatheta}
	\begin{split}
		&\int\frac{u^r}{r}(\omega^\theta)^2rdrdz\leq \|\frac{u^r}{r}\|_{L^\infty}\int(\omega^\theta)^2rdrdz=\epsilon\|\phi_{1,z}\|_{L^\infty}\|\omega^\theta\|_{L^2}^2\\
		&\leq C(\|\phi_{1,z}\|_{L^2}+\|\nabla^2\phi_{1,z}\|_{L^2})\|\omega^\theta\|_{L^2}^2\leq C(\|\phi_{1,z}\|_{L^2}+\|\nabla\omega_1\|_{L^2})\|\omega^\theta\|_{L^2}^2,
	\end{split}
\end{equation}
where we have used the Biot-Savart law \eqref{hjl} $\epsilon \phi_{1,z}(r,z)=\frac{u^r}{r}$, the Sobolev embedding, and the elliptic estimate \eqref{ellipticestimate} in the last step.

For the second nonlinear term in \eqref{eqwtheta}, we have 
\begin{equation} 
	\begin{split}
	\label{utheta4}
	& \left|\int\frac{(u^\theta)^2_z\omega^\theta}{r}rdrdz\right|=\left|\int \frac{(u^\theta)^2}{r}\omega^\theta_zrdrdz\right|\\
	& \leq \frac{1}{2}\int \frac{(u^\theta)^4}{r^2}rdrdz+\frac{1}{2}\|\nabla\omega^\theta\|^2.
\end{split}
\end{equation}
The first integral term in \eqref{utheta4} is estimated as 
\begin{equation}
	\label{endsep}
	\begin{split} \int\frac{(u^\theta)^4}{r^2}rdrdz&=\int_{r\geq 1}\frac{(u^\theta)^4}{r^2}rdrdz+\int_{r\leq 1}\frac{(u^\theta)^4}{r^2}rdrdz\\
		&\leq \int (u^\theta)^4rdrdz+{\int_{r\leq 1} u_1^4r^2rdrdz}\\
		&\leq \int (u^\theta)^4rdrdz+{\int_{r\leq 1} |u_1|^{4-\epsilon'}{|u_1|^{\epsilon'}}r^2rdrdz}\\
		&\leq C+C\|\Gamma^\epsilon_0\|_{L^\infty}^{\epsilon'}\leq C.
	\end{split}
\end{equation}
		
Adding up the estimates \eqref{firstomegatheta}, \eqref{utheta4} and \eqref{endsep} in \eqref{eqwtheta}, and using the Gronwell's inequality, we get that 
		\[ \|\omega^\theta(t)\|_{L^2}\leq C.\]
		Then since 
		\[ u^r\mathbf{e_r}+u^z\mathbf{e_z}=\epsilon \nabla \times (-\Delta)^{-1}(\omega^\theta \mathbf{e_\theta}),\]
		using Sobolev embedding, we have 
		\[ \|u^r\mathbf{e_r}+u^z\mathbf{e_z}\|_{L^6}\leq \|\nabla (u^r\mathbf{e_r}+u^z\mathbf{e_z})\|_{L^2}\leq C\|\omega^\theta\|_{L^2}\leq C.\]
		Then based on the {\em a priori} estimate \eqref{sae}, we have 
		\[ \|u^r\mathbf{e_r}+u^z\mathbf{e_z}\|_{L^4}\leq \|u^r\mathbf{e_r}+u^z\mathbf{e_z}\|^{\frac{3}{4}}_{L^6}\|u^r\mathbf{e_r}+u^z\mathbf{e_z}\|^{\frac{1}{4}}_{L^2}\leq C.\]
		
		This together with the estimate \eqref{utheta4e} and the Prodi-Serrin criterion \eqref{prodiserrin} implies the global regularity of the solutions. 

\end{proof}
\begin{remark}
	We proved our main result Theorem \ref{mainresult} using the $L^2$ estimate for $\omega_1$ and the $L^{4-\epsilon'}$ estimate for $u_1$. And we can also use the $L^p$ estimate for $\omega_1$ and $L^q$ estimate for $u_1$ with $q=2p-\frac{p\epsilon'}{2}$ for $p>1$ to get the same result.
\end{remark}

{\bf Acknowledgments.} The research was in part supported by the NSF Grants No. DMS-1613861 and DMS-1318377.
\bibliographystyle{plain}
\bibliography{selfsimilar}
\end{document}